\documentclass{amsart}
\usepackage{amsthm,amsmath,amssymb,amscd}
\newcommand{\thisdate}{\today}
{
   \newtheorem{theorem}{Theorem}[subsection]                     
   \newtheorem{proposition}[theorem]{Proposition}     
   \newtheorem{lemma}[theorem]{Lemma}

}
{\theoremstyle{definition}

   \newtheorem{definition}[theorem]{Definition}
   
}

\newcommand{\QQ}{{\mathbb{Q}}}

\newcommand{\PP}{{\mathbb{P}}}
\newcommand{\ZZ}{{\mathbb{Z}}}

\newcommand{\GG}{{\mathbb{G}}}

\newcommand{\bB}{{\mathbf{B}}}

\newcommand{\bK}{{\mathbf{K}}}
\newcommand{\bM}{{\mathbf{M}}}
\newcommand{\bmu}{{\boldsymbol{\mu}}}

\newcommand{\cB}{{\mathcal B}}
\newcommand{\cC}{{\mathcal C}}
\newcommand{\cD}{{\mathcal D}}

\newcommand{\cG}{{\mathcal G}}

\newcommand{\cK}{{\mathcal K}}
\newcommand{\cL}{{\mathcal L}}
\newcommand{\cM}{{\mathcal M}}

\newcommand{\cO}{{\mathcal O}}
\newcommand{\cP}{{\mathcal P}}

\newcommand{\cS}{{\mathcal S}}

\newcommand{\cX}{{\mathcal X}}

\newcommand{\Spec}{\operatorname{Spec}}
\newcommand{\Isom}{\operatorname{Isom}}

\newcommand{\Stab}{\operatorname{Stab}}

\newcommand{\Aut}{{\operatorname{Aut}}}

\newcommand{\das}{\dashrightarrow}
\newcommand{\dar}{\downarrow}
\newcommand{\ocL}{\overline{{\mathcal L}}}
\newcommand{\ocM}{\overline{{\mathcal M}}}
\newcommand{\obM}{\overline{{\mathbf M}}}

\newcommand{\sm}{{\operatorname{sm}}}
\newcommand{\bal}{{\operatorname{bal}}}
\newcommand{\adm}{{\operatorname{adm}}}
\newcommand{\rig}{{\operatorname{rig}}}
\newcommand{\double}{\genfrac..{0pt}1
{\raise -1pt\hbox{$\scriptstyle\to$}}{\raise 3pt\hbox {$\scriptstyle\to$}}}
\newcommand{\setmin}{\,\protect%
\begin{picture}(8,10)\qbezier(1,5.5)(4,4.)(7,2.5)\end{picture}\,}

\begin{document}
\title[Complete moduli for families over curves]{Complete moduli for families
   over semistable curves}
\author[D. Abramovich]{Dan Abramovich}
\thanks{D.A. Partially supported by NSF grant DMS-9700520 and by an Alfred
	P. Sloan research fellowship}  
\address{Department of Mathematics\\ Boston University\\ 111 Cummington
	 Street\\ Boston, MA 02215\\ USA} 
\email{abrmovic@math.bu.edu}
\author[A. Vistoli]{Angelo Vistoli}
\thanks{A.V.  partially supported by the University of
	Bologna, funds for selected research topics.}
\address{Dipartimento di Matematica\\ Universit\`a di Bologna\\Piazza di Porta
	 San Donato 5\\ 40127 Bologna\\ Italy}
\email{vistoli@dm.unibo.it}
\date{\thisdate}

\maketitle

This note is a research announcement, summarizing and explaining results proven
and detailed in forthcoming papers \cite{A-V:stable-maps},
\cite{A-V:fibered-surfaces}, \cite{A-C-J-V:covers}.
\section{Moduli and stacks} 
\subsection{Stacks as moduli objects} In the last two decades, it has been
observed that typically a ``nice'' moduli problem corresponds to a
Deligne-Mumford algebraic stack admitting a projective coarse moduli scheme.  
There are numerous examples of this phenomenon. Let us mention just a few:
\begin{enumerate}
\item $\ocM_{g,n}$: the moduli of stable $n$-pointed curves of genus $g$
\cite{Deligne-Mumford}; 
which, away from small characteristics generalizes to 
\item $\ocM_{g,n}(X,d)$: Kontsevich's moduli of stable $n$-pointed maps of
genus $g$ and degree $d$; and
\item $\cB G$: the moduli of principal homogeneous $G$-spaces, for a finite
group $G$. 
\end{enumerate}

Deligne-Mumford stacks form a 2-category which is an extension of the category
of schemes in a natural way. There is an extensive theory of cohomology,
intersection theory, vector bundles, and K-theory of stacks, which was
developed largely due to the importance of stacks in moduli theory. A closely
related notion of $\QQ$-varieties has been studied extensively. Some
natural diophantine equations related to stacks were studied in
\cite{Darmon-Granville}. 

\subsection{Stacks are basic objects} In this note we hope to
convince the 
reader that Deligne-Mumford stacks should be considered as basic objects of
algebraic geometry, like schemes, and not just as objects dedicated to
moduli problems. We argue as follows: a natural moduli problem of certain
stable families  over nodal curves is introduced; this moduli problem is not
complete; 
a natural compactification of this moduli problem involves families 
over {\em curves with Deligne-Mumford stack structure}; the resulting complete
moduli problem is a nice one, namely it is a complete Deligne-Mumford stack
admitting a projective coarse moduli scheme.

Let us introduce such moduli problems.

\subsection{The problem of moduli of families} Consider a Deligne-Mumford stack
$\cM$ 
admitting a projective coarse moduli scheme $\bM\subset \PP^N$. Given a curve
$C$, it is often natural to consider morphisms $f:C \to \cM$ (or equivalently,
objects $f\in \cM(C)$): in case $\cM$ is the moduli of geometric objects, these
morphisms correspond to families over $C$. For example, if
$\cM=\ocM_{\gamma}$, then 
morphism $f:C\to \cM$ correspond to families of stable curves of genus
$\gamma$ over $C$; and if 
$\cM = \cB G$ we get principal $G$-bundles over $C$. It should be obvious that
it is interesting to study moduli of such objects; moreover, it is natural
to study such moduli as $C$ varies, and find a natural compactification for
such moduli. 

\subsection{Stable maps} Denote by $g$ the genus of $C$.
In case $\cM$ is represented by a projective scheme $X\subset \PP^N$, a natural
answer to these questions is given by the Kontsevich stacks of stable maps
$\ocM_g(X, d)$. It is tempting to mimic this construction in the case of an
arbitrary stack as 
follows: let $C$ be a nodal projective connected curve; then a morphism $C \to
\cM$  is said to be a {\em stable map} of degree $d$ if the associated morphism
to the coarse moduli scheme $C \to \bM$ is a stable map of degree $d$.  

It follows from our results below that this moduli problem is  a
Deligne-Mumford stack. A somewhat surprising point is, that it is not
complete.

To see this, we fix $g=2$ and consider the specific case of $\cB G$ with $G=
(\ZZ/3\ZZ)^4$. Any 
smooth curve of genus $2$ admits a {connected} principal $G$ bundle,
corresponding of a surjection $H_1(C, \ZZ) \to G$. If we let $C$ degenerate to
a nodal curve $C_0$ of geometric genus $1$, then   $H_1(C_0, \ZZ) \simeq
\ZZ^3$, and since there is no surjection $\ZZ^3 \to G$, there is no connected
principal $G$-bundle over $C_0$. This means that there is no limiting stable
map  $C_0 \to \cB G$.

\subsection{Enter orbispace-curves}

Our main goal here is to correct this deficiency. In
order to do  so, we will enlarge the category of stable maps
into $\cM$. The source curve
${\cC}$ of a new stable map ${\cC} \to \cM$ will
acquire an orbispace structure at its nodes. Specifically, we
endow it  with the structure of a Deligne-Mumford stack.

It is not hard to see how these orbispace structure come
about. Let
$S$ be the spectrum of a discrete valuation ring $R$ of pure characteristic 0,
with quotient field $K$, and let
$C_K\to \eta\in S$ be a nodal curve over the generic point, together
with a map $C_K \to \cM$ of degree $d$, whose associated map $C_K \to \bM$ is
stable. We can 
exploit the fact that
${\ocM}_{g,0}(\bM,d)$ is complete; after a ramified base
change on
$S$ the induced map $C_K \to
\bM$ will extend to a stable map $C \to \bM$ ove $S$. Let
$C_\sm$ be the smooth locus of the morphism $C \to S$;
Abhyankar's lemma, plus a fundamental purity lemma (see
\ref{Lem:purity} below) shows that after a suitable base change we can
extend the map $C_K \to
\cM$ to a map $C_\sm \to \cM$; in fact the purity lemma fails
to apply only at the ``new'' nodes of the central fiber, namely those which are
not in the closure of nodes in the generic fiber. On the other hand, if
$p\in C$ is such a node, then on an
\'etale neighborhood $U$ of
$p$, the curve $C$ looks like $$uv = t^r,$$ where $t$ is the
parameter on the base. By taking
$r$-th roots:
$$u = u_1^r;\, v = v_1^r$$
we have a {\em nonsingular} cover
$V_0\to U$ where $V_0$ is defined by
$u_1v_1 = t$. The purity lemma applies to $V_0$,
so the composition ${V_0}_K \to C_K \to \cM$ extends over all
of $V_0$. There is a minimal intermediate cover $V_0\to
V\to U$ such that the family extends already over
$V$; this $V$ will be of the form $xy = t^{r/m}$, and the
map $V \to U$ is given by $u = x^m$,
$v = y^m$. Furthermore, there is an action of the group
$\bmu_m$ of roots of 1, under which $\alpha\in \bmu_m$ sends $x$
to $\alpha x$ and $y$ to $\alpha^{-1} y$, and
$V/ \bmu_m = U$. This gives the orbispace structure ${\cC}$
over $C$, and the map $C_K \to
\cM$ extends to a map ${\cC} \to \cM$.

This gives the flavor of our definition, which we will  give below in a general
setting. 

Here is the lemma we used in the argument:

\begin{lemma}[Purity Lemma, \cite{A-V:fibered-surfaces}]\label{Lem:purity} Let $\cM$ be a
separated Deligne-Mumford stack, $\cM \to
\bM$ its coarse moduli space. Let $X$ be a separated
scheme of dimension 2 satisfying Serre's condition {\rm S2}.
Let
$P\subset X$ be a finite subset consisting of closed points,
$U=X\setmin P$, and assume that the local fundamental groups
of
$U$ around the points of $P$ are trivial. Then a morphism
$U \to \cM$ extends to a morphism $X \to \cM$ if and only if the
composition $U \to \cM \to \bM$ extends to a morphism $X \to
\bM$.
\end{lemma}

\subsection{Restrictions on the residue characteristics of the base scheme} As
indicated above, we 
will need to apply Abhyankar's lemma. This means 
that we need to assume that no characteristic $p$
appearing divides the order of a stabilizer of a geometric point  of $\cM$ of
characteristic $p$. This requirement is enough to guarantee that the
moduli category described below
is an Artin stack with finite diagonal. To get a Deligne-Mumford stack we need
a bit more, to ensure that 
$\ocM_{g,n}(\bM,d)$ is a Deligne-Mumford stack: a-priori it is an Artin stack
with finite diagonal, but there is
a dense open set of primes in $\ZZ$, depending on $g,n,d$ and $\bM$,  over
which $\ocM_{g,n}(\bM,d)$ is a  Deligne-Mumford stack. 
 If we denote by $e(M)$ the product of all the ``bad'' primes listed above,
we require that all schemes considered below are schemes over
$\Spec\ZZ[1/e(M)]$. 
(For simplicity the reader may wish to
stick with a base scheme of characteristic 0.)

\subsection{Twisted unpointed nodal curves} Our first goal is to identify what
type of ``orbispace curves'' we want to work with. It is convenient to start
with ``unpointed'' curves. 

\begin{definition}
A {\em twisted nodal curve over $S$} is a diagram
$$\begin{array}{c} \cC  \\ \dar  \\ C  
\\ \dar  \\ S  
\end{array}$$
Where
\begin{enumerate}
\item $\cC$ is a Deligne-Mumford stack, with geometrically connected fibers
over $S$, which is \'etale locally a nodal curve
over $S$;
\item the morphism $\cC \to C$ exhibit $C$ as the coarse moduli scheme of
$\cC$; and
\item $\cC \to C$ is isomorphic away from the nodes.
\end{enumerate}
\end{definition}

In other words, a twisted nodal curve is a way to put a Deligne-Mumford
stack structure $\cC$ ``over'' a nodal curve $C$. One can give an explicit
description of such a structure. For instance, if $S$ is the spectrum of an
algebraically closed field $k$, then locally in the \'etale topology $\cC$ is the
stack quotient of $\Spec k[x,y]/(xy)$ by the action of $\bmu_m$, where $m$ is
prime 
to the characteristic of $k$, the parameters $x$ and $y$ are eigenvectors of
the action, 
and the eigenvalues are {\em primitive} $m$-th roots of $1$.

\subsection{Twisted pointed nodal curves}
Recall that a natural way to get a {\em pointed}  nodal curve $C$ from an
unpointed one, 
is obtained by ``separating'' some of the nodes and then ``ordering'' the
points 
above these nodes. These points are disjoint  sections
of the smooth locus of the curve.  

If one ``separates'' a node on a {\em twisted} pointed curve, one obtains an
object which is a bit more subtle than a section. To see what happens we can
look at an \'etale neighborhood of such a ``separated node'', which is the
stack quotient of ${\cS}pec_S\cO_s[x]$ by a faithful action of $\bmu_m$ on the
variable 
$x$. The quotient of the section $\{x=0\}$ in this \'etale neighborhood is a
copy of the classifying stack $\cB\bmu_m/S$. When these \'etale neighborhoods
are glued together we obtain an {\em \'etale gerbe over 
$S$.} For our purposes, the reader may think of an \'etale gerbe  $\cG \to S$
as a stack, which locally 
in the \'etale topology is isomorphic to the classifying stack $\cB G/S$ for
some finite \'etale group scheme $G \to S$. (Formally, $\cG \to S$ is an
\'etale gerbe when both $\cG \to S$ and the diagonal $\cG \to\cG\times_S\cG$
are surjective \'etale morphisms.) 
Thus, locally in the \'etale topology 
there is a section $S \to \cG$, but such sections may not exist globally.

This motivates the following definition:

\begin{definition}
A {\em twisted nodal $n$-pointed curve over $S$} is a  diagram
$$\begin{array}{ccc} \cC_i^s & \subset & \cC  \\  &\searrow & \dar  \\ &
 & C   
  \\ &  & \dar  \\ &&S  
\end{array}$$
Where
\begin{enumerate}
\item $\cC$ is a Deligne-Mumford stack which is \'etale locally a nodal curve
over $S$;
\item $\cC_i^s  \subset  \cC$ are disjoint closed substacks in the smooth locus
of $\cC \to S$;  
\item $\cC_i^s \to S$ are \'etale gerbes; 
\item the morphism $\cC \to C$ exhibit $C$ as the coarse moduli scheme of 
$\cC$; and
\item $\cC \to C$ is isomorphic away from the nodes and the  $C_i^s$.
\end{enumerate}
\end{definition}

Note that if we let $C_i^s$ be the coarse moduli spaces of $\cC_i^s$, then
$C_i^s$ embed in $C$ - they are the images of $\cC_i^s$, and $C$ becomes a
usual nodal pointed curve. We say that $\cC \to S$ is a twisted pointed curve
of genus $g$, if $C \to 
S$ is a pointed curve of genus $g$.

\subsection{Morphisms of twisted pointed nodal curves}

\begin{definition} Let $\cC \to S$ and $\cC' \to S'$ be twisted $n$-pointed
curves. A morphism $F : \cC \to \cC'$ is a cartesian diagram
$$\begin{array}{ccc}
\cC  & \stackrel{F}{\to} & \cC' \\
\dar &                   & \dar \\
S    & \stackrel{f}{\to} & S' 
\end{array}$$
such that $F^{-1}{\cC'}_i^s = \cC_i^s$.
\end{definition}

Since twisted pointed curves are stacks rather than schemes, we need to be a
bit careful. 
If $F,F_1 : \cC \to \cC'$ are morphisms, then we can define a 2-morphism $F \to
F_1$ to be an isomorphism of functors. In this way, twisted pointed curves form
a 2-category.  This may seem to be
a problem, since we wish to use them to form a stack, which, by definition, is
a category. Here the fact that $\cC \to C$ is generically isomorphic comes to 
rescue: it is easy to see that 2-morphisms are unique when they exist, and
replacing morphisms by their equivalence classes we have the following:

\begin{proposition}[\cite{A-V:stable-maps}] The 2-category of twisted pointed
curves is equivalent to a category. 
\end{proposition}
We call the resulting category {\em the category of twisted pointed curves}.

\subsection{Stable maps into a stack} As before, we consider a proper
Deligne-Mumford stack $\cM$ admitting a projective coarse moduli scheme
$\bM$. We fix a projective embedding $\bM\subset \PP^N$.

\begin{definition}
A {\em twisted stable $n$-pointed map of genus $g$ and degree $d$ over $S$}
$$(\cC \to S, \cC_i^s\subset \cC, f: \cC \to \cM)$$  
consists of a commutative diagram
$$\begin{array}{ccc} \cC &\to& \cM \\
		\dar & & \dar \\
		C & \to & \bM \\
\dar &&\\
S &&
\end{array}$$
 along with $n$ closed substacks $\cC_i^s\subset \cC$,
satisfying:
\begin{enumerate}
\item $\cC \to C \to S$ along with $\cC_i^s$  is a twisted nodal $n$-pointed
curve 
over $S$; 
\item the morphism $\cC \to \cM$ is representable; and
\item $(C\to S, C_i^s, f:C \to \bM)$ is a stable $n$-pointed map of
degree $d$.   
\end{enumerate}
\end{definition}

A few remarks are in order.
\begin{enumerate}
\item the prefix ``twisted'' comes to stress the fact that the base curve $\cC$
has ``extra structure'' as a  Deligne-Mumford stack. A twisted stable map where
$\cC \to C$ is an isomorphism is called ``untwisted''.
\item Twisted stable maps $\cC \to \cM$ can be defined without invoking the
coarse moduli scheme $C$. For instance, the stability of $C \to \bM$ is
equivalent 
to the assertion that $\Aut_{\cM}(\cC, \cC_i^s)$ is finite.
\item The condition that the morphism $\cC \to \cM$ be representable means that
the stack structure on $\cC$ is the minimal necessary to ensure the existence
of the morphism $\cC \to \cM$.  This should be considered a stability
condition, in the sense that it is essential to ensure that the moduli
problem be separated.  
\end{enumerate}

Now we define morphisms of twisted stable maps:   

\begin{definition} A {\em morphism of twisted stable maps} 
$$ G: (\cC \to S, \cC_i^s, f: \cC \to \cM) \to (\cC' \to S',
 {\cC'}_i^s, f': \cC' \to \cM)$$
 consists of data $G=(F,\alpha)$, where $F:\cC \to \cC'$ is a morphism of
twisted pointed curves, and  $\alpha: f \to f'\circ F$ is an isomorphism. 
\end{definition}

Note that, unlike stable maps into a scheme, a twisted stable map $f:\cC \to
\cM$ may have automorphisms which are trivial on the source $\cC$, even when
$\cC = C$. For example, a family of stable curves over $C$ may have
automorphisms fixing $C$. This is the role of $\alpha$ in the definition.

Again, twisted stable maps naturally form a 2-category. But by the proposition,
this 2-category is equivalent to a category.
We call this category {\em the category of twisted stable maps}. In
\cite{A-V:stable-maps} we also give an explicit realization
of this category, which is unfortunately a bit technical, in terms of atlases
of charts over the coarse curves $C$. This is in 
analogy with Mumford's treatment of $\QQ$-varieties in
\cite{Mumford:towards}. Both descriptions are useful when proving results
about this category.  

It is natural to denote this category $\ocM_{g,n}(\cM,d)$, but we find the
abundance of $\cM$'s a bit confusing. We propose to denote it instead by
$\cK_{g,n}(\cM,d)$.  

\subsection{The main result} There is structural functor
$\cK_{g,n}(\cM,d) \to \mathcal{S}ch$ which
associates to a twisted stable map $ (\cC \to S, \cC_i^s, f: \cC \to
\cM)$ 
the base scheme $S$. With this functor, our main result is:
\begin{theorem}
The category $\cK_{g,n}(\cM,d)$  forms a Deligne-Mumford stack, admitting a
projective coarse moduli scheme $\bK_{g,n}(\cM,d)$.
\end{theorem}

The proof of the theorem is far from easy. After all, many classical moduli
problems are solved by stack quotients of appropriate Hilbert schemes, and 
analogues for Hilbert schemes for stacks are not simple to construct. Our
construction builds on the 
fact that the Kontsevich stack $\cK_{g,n}(\bM,d)$ is known to be a complete
Deligne-Mumford stack, with projective moduli space.

The main steps in our proof are the following.
\begin{enumerate}
\item  We prove that the diagonal of $\cK_{g,n}(\cM,d)$ is finite and
representable.

\item  We show that $\cK_{g,n}(\cM,d)$ is an algebraic stack, by
checking that the condition of \cite{Artin:stacks} are verified. Here the
hard part is the proof of algebraizability of formal twisted stable
maps; this involves a form of Grothendieck's existence
theorem for algebraic stacks.

\item  We check the valuative criterion for properness for
$\cK_{g,n}(\cM,d)$. For this the main tool is Abhyankar's lemma,
together with the purity lemma.

\item  We prove boundedness; this is based on showing that $\cK_{g,n}(\cM,d)
\to \cK_{g,n}(\bM,d)$ has finite fibers, but is a bit more involved. Together
with the valuative criterion this implies properness, in particular
$\bK_{g,n}(\cM,d) \to \bK_{g,n}(\bM,d)$ is finite, thus $\bK_{g,n}(\cM,d)$ is
projective. 
\end{enumerate}

\subsection{Balanced maps} When we first introduced orbispace curves into the
picture, the extra structure appeared at a ``new'' node  in the central fiber
only in the following way: near a node locally of the form $U: uv=t^r$ we had
an orbispace chart, locally of the form  $V: xy = t^{r/m}$ with a $\bmu_m$
action such that $U = V/\bmu_m$. The action was given as follows: $\alpha\in
\bmu_m$ sends $(x,y)\mapsto (\alpha x , \alpha^{-1} y)$. Note that this is
not the most general action of the stabilizer of a node appearing in a twisted
stable map: here the eigenvalues of the stabilizer acting on the tangent spaces
to the two
branches of the node are inverse to each other. A
twisted stable map with this property is
called {\em balanced}. We denote the subcategory of balanced twisted stable
maps $\cK_{g,n}^{\bal}(\cM,d)$. We have:

\begin{proposition}
The subcategory $\cK_{g,n}^{\bal}(\cM,d)\subset \cK_{g,n}(\cM,d)$ is an open
and 
closed substack. It contains the closure of the locus of twisted stable
pointed maps
with nonsingular source curve $\cC$.
\end{proposition}

\section{Fibered surfaces}
\subsection{Fibered surfaces and coarse fibered surfaces} As our first example,
we consider the case $\cM = \ocM_{\gamma,\nu}$ of stable $\nu$-pointed curves
of genus $\gamma$. For simplicity we look at the case $n=0$ and omit $n$ from
the notation.  

An untwisted stable map $C \to \ocM_{\gamma,\nu}$ is equivalent to a family $(X
\to C,\tau_i:C \to X)$ of stable $\nu$-pointed curves of genus $\gamma$ over a
nodal curve $C$, 
such that $\Aut(X\to C, \tau_i)$ is finite. In other words, $X$ is a surface,
mapping to $C$ with sections, with  stable pointed fibers, satisfying a further
stability  
condition. In general, a twisted stable map $\cC \to \ocM_{\gamma,\nu}$ is
equivalent to a similar family $\cX \to \cC$, only now $\cX$ and $\cC$ are not
schemes but 
stacks. Nevertheless, we call the family  $\cX \to \cC$, associated to a
twisted stable map, a {\em fibered surface}. Thus fibered surfaces appear
naturally in the boundary of the moduli of untwisted fibered surfaces.

Given a fibered surface  $\cX \to \cC$ one may consider the coarse moduli
schemes:

$$\begin{array}{ccc} \cX &\to &\cC\\
		\dar & & \dar \\
		 X & \to & C.
\end{array}$$

Now $X \to C$ is not necessarily a family of stable pointed curves; rather, it
is locally the {\em quotient} of such a family by the action of a cyclic
group. We call it a {\em coarse fibered surface}.

\subsection{Comparison with Alexeev's work} One may ask, what are the
singularites of coarse fibered surfaces? And, to what 
extent can coarse fibered surfaces replace fibered 
surfaces in the boundary of moduli? In some sense, these questions have already
been addressed in the literature. If one restricts attention to
{\em balanced} coarse fibered surfaces, then we show:

\begin{proposition}[\cite{A-V:fibered-surfaces}]
\begin{enumerate}
\item  Let  $\cX \to \cC$ be a balanced fibered surface, $X \to C$ the
associated  coarse fibered surface. Then $X$ has semi-log-canonical
singularities. 
\item Consider the morphism $X \to \obM_{\gamma,\nu}$ obtained by composing $X
\to C$ with the structural morphism $C \to \obM_{\gamma,\nu}$. This morphism is
a stable map in the sense of Alexeev \cite{Alexeev:stable-maps}.
\item In characteristic 0, there is a finite morphism from
$\cK_{g}^{\bal}(\ocM_{\gamma,\nu},d)$ to Alexeev's moduli stack of stable
maps. 
\end{enumerate}
\end{proposition}

One can use this to define {\em the  stack of coarse fibered surfaces} as
follows. A coarse fibered surface $X \to C \to S$ over a base scheme $S$ gives
rise to an Alexeev stable map $X \to \obM_{\gamma,\nu}$ as well as a Kontsevich
stable map $C \to \obM_{\gamma,\nu}$; and in addition we have the  morphism $X
\to C$.  Since both the Alexeev stacks and Kontsevich stacks are
Deligne-Mumford stacks, it is an easy consequence of Grothendieck's theory of
the Hilbert scheme that there is a Deligne-Mumford stack $A$ of such triples
$$(X \to \obM_{\gamma,\nu},C \to \obM_{\gamma,\nu},X \to C),$$ admitting a 
quasi-projective coarse moduli scheme. There is a morphism 
$\cK_{g}^{\bal}(\ocM_{\gamma,\nu},d) \to A$, and the image is the  stack of
coarse fibered surfaces. 

We also show that for suitable choices of the parameters, this morphism is not
one-to-one, and is ramified. This means that some of the singularities of
Alexeev-type stacks are partially resolved in our stack. In view of this, 
one can say that the stack of balanced fibered surfaces is, in a sense, a
refinement of Alexeev's work on surface stable maps, in the particular case
described here.

\subsection{Towards plurifibered varieties} Our main theorem has a nice
recursive feature: the input is a Deligne-Mumford stack with a projective
coarse moduli scheme, and the output is of the same nature. It is tempting to
apply this feature to higher dimensional varieties. Given a sequence of
dominant 
rational maps $$X_n \das X_{n-1} \das \cdots \das X_1 \das S$$ of relative
dimension 1,
one can
apply our construction inductively and get, first, a {\em canonical model}
$$\overline{X_n}_\eta \to \overline{X_{n-1}}_\eta \to \cdots \to
\overline{X_1}_\eta  \to \eta$$ over the generic point $\eta\in S$, and,
after replacing $S$ by a suitable proper, generically finite and surjective
base change, one can get a {\em stable reduction} 
$$\overline{X_n} \to \overline{X_{n-1}}\to \cdots \to
\overline{X_1} \to S.$$ We call these structures {\em stable plurifibered
varieties}. 
Alexeev has suggested to use the minimal model
conjecture to define stable varieties (as well as canonical models) in a
general situation, without the
presence of a ``plurifibration'' as above. Interestingly, the structure
given by the plurifibration allows one to
bypass the minimal model program entirely. It would be interesting to compare
the singularities of stable plurifibered varieties to the singularities which
arise in the minimal model program.

\section{Twisted covers and level structures}

\subsection{Twisted principal bundles}

We fix a finite group, or a finite \'etale group-scheme $G$, and we set $\cM =
{\cB}G$. We denote the stack ${{\cK}_{g,n}({\cB}G,0)}$ by
$\cB_{g,n}G$, and ${{\cK}_{g,n}^{\bal}({\cB}G,0)}$
by $\cB_{g,n}^{\bal}G$. These we call the stack of {\em stable twisted
$G$-bundles}  and {\em balanced stable twisted $G$-bundles}. These names are
motivated by the following observation:

There is a nice explicit description of twisted stable maps to $\cB G$ in
terms of schemes. Say we are working over a base scheme $T$ and $G \to T$
is a finite \'etale group scheme. Then the associated morphism $T \to \cB G$
exhibits $T$ as the universal principal $G$ bundle $\cP G \to \cB G$. 
Let $\cC \to \cB G$ be a twisted stable map. Pulling back the universal
principal bundle $\cP G \to \cB G$ we get an associated {\em twisted
principal $G$-bundle} $P \to \cC$. The fact that $\cC \to \cB G$ is
representable and $\cP G$ is a scheme means that $P$ is a scheme. Since $P\to
\cC$ is \'etale and $\cC$ is nodal, $P$ is a nodal curve, not necessarily
connected. Moreover, the 
action of $G$ on $P$ has no fixed points away from the nodes and the points
lying above the marked points of $\cC$, and the schematic quotient is $C$. On
the other hand,
given such a Galois cover $P \to C$, we can recover $\cC$ as the stack quotient
of $P$ by $G$. One easily checks the following:

\begin{proposition}
The category $\cB_{g,n}G$ is equivalent to the category of stable $n$-pointed
curves $C$ of genus $g$, along with nodal $G$-covers $P \to C$ such that the
action is fixed-point free away from the nodes and the marked points.
\end{proposition}

The local structure of these stacks is simple.
A straightforward calculation of deformation and obstruction spaces, similar to
the one in \cite{Deligne-Mumford}, gives

\begin{theorem}
The stack $\cB_{g,n}G$ is smooth over $\ZZ[1/|G|]$. Its relative dimension at a
given twisted stable map $\cC \to \cB G$  is $3g-3+n-r$, where $r$ is the
number of nodal 
points at which $\cC$ is not balanced.
\end{theorem}

\subsection{Balanced twisted covers and Galois admissible covers} In
\cite{Harris-Mumford}, Harris and Mumford considered a compactification of the
Hurwitz space of simply branched covers of $\PP^1$ via {\em admissible
covers}. The construction generalizes to covers of curves of arbitrary genus
$g$ with arbitrary ramification type, see  \cite{Mochizuki:Hurwitz} and
\cite{Wewers:Hurwitz}. We will now 
relate these generalized Harris-Mumford stacks of admissible covers with our
construction.

The stack
$\cB_{g,n}^{\bal} G$ contains an open and closed substack parametrizing
connected $G$ bundles. The condition for the Galois cover  $P \to C$ to be
balanced is equivalent to the Harris-Mumford condition for an admissible
cover. We have the following: 

\begin{proposition} The subcategory $\cB_{g,n}^{\adm} G\subset \cB_{g,n}^{\bal}
G$ of connected balanced principal bundles is an open and closed substack,
which is isomorphic to the stack of Galois admissible covers with Galois group
$G$. 
\end{proposition}

Note that while the stack described by Harris-Mumford is in general singular,
the proposition implies that the stack of Galois admissible covers is
smooth. This fact was also observed by Wewers
\cite{Wewers:Hurwitz}. Remarkably, while Wewers's approach differs from ours,
the resulting stacks are the same. 

\subsection{Twisted \'etale covers and admissible covers} In order to treat
admissible covers which are not Galois, we can use the equivalence of
categories 
between $d$-sheeted \'etale covers and principal $\cS_d$-bundles, where $\cS_d$
is the symmetric group on the $d$ letters $\{1,\ldots,d\}$.  Given a branched
cover $D \to C$ of a smooth base curve of genus $g$, marked by the $n$ branch
points, we take the associated Galois cover $P \to C$ with Galois group
$\cS_d$. This is an object of $\cB^{\bal}_{g,n}\cS_d$. Thus, a natural
compactification of the moduli of such branched covers is an open and closed
substack $Adm_{g,n,d}$ in $\cB^{\bal}_{g,n}\cS_d$. 

\begin{proposition} Given an object $P \to C$ of
$\cB^{\bal}_{g,n}\cS_d$, the schematic quotient $D = \cS_{d-1} \backslash P $
is a 
$d$-sheeted cover of $C$, and $P \to C$ is in $Adm_{g,n,d}$ if and only if $D$
is 
connected. The branched cover $D \to C$ is admissible in the sense of
Harris-Mumford. The stack $Adm_{g,n,d}$ is the normalization of the stack of
generalized Harris-Mumford admissible covers.
\end{proposition}

Since the stack of Harris-Mumford admissible covers is singular, it is not
isomorphic to $Adm_{g,n,d}$. In other words, $Adm_{g,n,d}$ is a minimal
desingularization of 
the Harris-Mumford stack.

There is a definition of $Adm_{g,n,d}$ which does not invoke the principal
bundle, only the covering $D \to C$ with extra structure. Consider the
following category of {\em twisted admissible covers:} objects over $S$ consist
of  
\begin{enumerate}
\item a balanced twisted curve $\cC \to S$; and
\item a connected finite \'etale cover  $\cD \to \cC$ of degree $n$,
\end{enumerate}
satisfying the follwing stability conditions: 
\begin{enumerate} 
\item the morphism $\cD \to \cC$ is representable;
\item the coarse curve $C \to S$ of $\cC$ is stable; and
\item for any geometric point $p$ of $\cC$, the action of $\Aut p$ on the fiber
$\cD_p$ is effective.
\end{enumerate}
We define morphisms by fibered diagrams as usual. We have:
\begin{proposition}
The category of twisted admissible covers  is a stack, isomorphic to
$Adm_{g,n,d}$. 
\end{proposition}

\subsection{Principal $\bmu_m$-bundles and invertible sheaves}

Stable twisted $\bmu_m$-bundles have a natural description via invertible
sheaves. We define a category as follows: the objects consist of data 
$(\cC, \cL, f)$,  where $\cC$ is a nodal twisted curve, $\cL \to \cC$ is an
invertible sheaf, and $f:\cL^m \to \cO_{\cC}$ is an isomorphism. We need a
stability condition: we require the coarse curve $C$ to be stable, and for each
node $p$ on $\cC$, we require the action of the stabilizer $\Stab_p$ on $\cL_p$
to be faithful. We call these objects {\em stable twisted $m$-torsion
invertible sheaves}. Morphisms of such objects are defined as fibered
diagrams. There is a  notion of 2-morphisms, making this into a 2-category, but
as before it is equivalent to a category. 

We claim that these are nothing but Stable twisted $\bmu_m$-bundles. Denote by
${\boldsymbol 1}$ the identity section in the total space of the bundle
$\cO$. Given a stable twisted $m$-torsion  invertible sheaf $(\cC, \cL, f)$, 
the inverse image $P=f^{-1}{\boldsymbol 1}$ in the total space of $\cL$ is a
principal $\bmu_m$-bundles,  and the stability condition on $(\cC, \cL, f)$
implies that  
$P \to \cC$ is stable. On the other hand, given a $\mu_m$ bundle one has an
associated $\GG_m$ bundle  extending to an invertible sheaf. It is easy to
verify that this is an equivalence of categories.

\subsection{Abelian level structures} There is a remarkable application of the
stacks $\cB_g G$ to Mumford's moduli of curves with level structures. In this
section, we assume that the structure sheaf of the base scheme contains 
the $m^{\rm th}$ roots of 1, and that we have fixed an
isomorphism $\bmu_m \simeq {\bf Z}/m$. We
will construct a smooth complete Deligne--Mumford stack
$\cL_g^m$, endowed with a finite morphism
$\cL_g^m \to \ocM_g$, which coincides with the
scheme of level $m$-structures over the stack ${\cM}_g$ of
smooth curves of genus $g$.

The idea is that we can interpret,
via Poincar\'e duality, a level $m$-structure on a smooth
curve $C$ as an element of the cohomology group $H^1\left(C,
(\ZZ/m)^{2g}\right)$; this in turn corresponds to an isomorphism class
of $(\ZZ/m)^{2g}$-bundles. This suggest that we can use our stack
${\cB}_{g,0}(\ZZ/m)^{2g}$ to define a level $m$-structure. 

Of course there is a problem here: a twisted bundle
with group  $(\ZZ/m)^{2g}$ over a fixed smooth curve has $(\ZZ/m)^{2g}$ as its
automorphism group. This is 
in contrast with the fact that the moduli stack of smooth curves with level
structure is representable. This is the same problem that one encounters  with
the Picard scheme: the stack of line bundles is not representable, because
every line bundle has $\GG_m$ as automorphisms, and one goes through a process
of  ``removing'' this $\GG_m$ action and sheafifying
(see \cite{FGA}, I.B.4, II.C.3, V.1).

This procedure can be carried out in general.
\begin{proposition} Let ${\mathcal
X}$ be a Deligne--Mumford stack. Suppose the automorphism group of every object
of ${\mathcal  X}$ contains a fixed  subgroup $G$, and that the embedding of
this subgroup commutes with base changes. Then 
there exists a Deligne--Mumford stack ${\mathcal Y}$, equipped with a
morphism ${\mathcal X} \to {\mathcal Y}$ which makes ${\mathcal X}$ into
an \'etale gerbe over ${\mathcal Y}$, so that the isomorphism
classes of geometric points are the same, but the
automorphism group of an object of ${\mathcal Y}$ is the
automorphism group of an object of ${\mathcal X}$, divided out by
$G$. 
\end{proposition}
This fact is certainly known, but we do not know a reference.
One can see this using \'etale presentations, as follows.
Take an \'etale map of finite type $U \to {\mathcal X}$, and set
$R = U \times_{\mathcal X}U$, so that $R \double U$ is an \'etale
presentation for ${\mathcal X}$. If $\xi\in{\mathcal X}(U)$ is the
object corresponding to the morphism $U \to {\mathcal X}$, then
$R$ represents the functor $\Isom_{U \times_S U}({\rm
pr}_1^*\xi, {\rm pr}_2^*\xi)$, where the ${\rm pr}_i \colon U
\times_S U \to U$ are the two projections. There is a free
action of $G$ on $R$, leaving the two projections $R \to U$
invariant, defined by composing isomorphisms with the 
automorphisms associated with $G$. This allows to define a 
quotient \'etale groupoid $R/G \double U$; this is an
\'etale presentation of the stack ${\mathcal Y}$. 

Let $G$ be  a finite abelian group, and consider $\cB_{g,n}^{\bal} G$. Every
object of this stack has $G$ in its automorphism group. Applying the
proposition, we obtain a stack $\cB_{g,n}^{\rig} G$, called the stack of
{\em rigidified} balanced twisted bundles, and a morphism
$\cB_{g,n}^{\bal} G\to \cB_{g,n}^{\rig} G$ as above. Given an {\em irreducible}
curve $C$ and a twisted $G$ bundle $P \to C$, the automorphism group of $P$
over $C$ is equal to $G$. This means that  the restriction $\rho\colon
\rho^{-1}\overline 
{\cM}'_{g,n}
\to \overline {\cM}'_{g,n}$ of the morphism $\rho\colon
\cB_{g,n}^{\rig}G \to \overline {\cM}_{g,n}$ to the open substack
$\overline {\cM}'_{g,n} \subseteq\overline {\cM}_{g,n}$
of {\em irreducible} stable curves is representable.

Consider the case $n=0$, $G = (\ZZ/m\ZZ)^{2g}$. In this case it is easy to
compare the stack $\cB_{g}^{\rig} G$ with Mumford's moduli scheme $\cM_g^{(m)}$
of curves with symplectic level-$m$ structure. 

\begin{proposition}
There is an open embedding $\cM_g^{(m)}\subset \cB_{g}^{\rig}
(\ZZ/m\ZZ)^{2g}$. The closure $\ocM_g^{(m)}$ of $\cM_g^{(m)}$ is an open and
closed substack, which coincides with the normalization of $\ocM_g$ in
$\cM_g^{(m)}$.  
\end{proposition}

The stack $\ocM_g^{(m)}$ can be described directly using symplectic
structures. The point is that, for a twisted curve $\cC$ underlying an object
in $\ocM_g^{(m)}$ over an algebraically closed field, we have
$H^1_{\mbox{\tiny{\'et}}}(\cC, 
\ZZ/m\ZZ) \simeq (\ZZ/m\ZZ)^{2g}$, and moreover this group carries a canonical
symplectic structure. We can also characterize the ``amount of twisting''
needed for such a curve: say $\cC$ is a {\em pre-level-$m$ balanced curve} if
the stabilized at each separating node is trivial and the stabilizer at a
non-separating node is cyclic of order $m$. We can define a category $\ocL_g^m$
of {\em twisted curves with level $m$ structure} whose objects are families
pre-level-$m$ 
balanced  
curves $\pi:\cC\to S$ along with  symplectic isomorphisms $R^1\pi_*\ZZ/m\ZZ
\to(\ZZ/m\ZZ)^{2g}_S$, and morphisms given by fibered squares as usual. This
forms a stack, and we have

\begin{proposition} The stack of twisted curves with level $m$  structure is
isomorphic to $\ocM_g^{(m)}$.
\end{proposition}

\subsection{Non-abelian  level structures} We return to an arbitrary finite
group $G$. Consider the stack $\cB_{g}^{\adm} G$ of connected balanced twisted
bundles. Objects $P \to C$ over smooth curves correspond to epimorphisms
$\pi_1(C) \to G$, once one chooses a base point on $C$. We can view the stack 
 $\cB_{g}^{\adm} G$  as a stack of stable twisted curves with Teichm\"uller
level structure with group $G$. In view of the work of Looijenga
\cite{Looijenga:level} and Pikaart-De Jong \cite{ Pikaart-DeJong:level}, it
would be interesting
to describe the coarse moduli scheme $\bB_{g}^{\adm} G$ and its relation with 
$\cB_{g}^{\adm} G$. This is the content of work in progress we are conducting
with A. Corti and J. de Jong.

\end{document}